\documentclass[11pt]{article}
\usepackage{eurosym}
\usepackage{amsfonts}
\usepackage{amsmath}
\usepackage{geometry}
\usepackage{amssymb}
\usepackage{graphicx}
\usepackage{color}
\usepackage[author={Oana}]{pdfcomment}
\usepackage[colorinlistoftodos]{todonotes}
\usepackage{titling}
\usepackage{lipsum}

\newtheorem{theorem}{Theorem}

\newtheorem{definition}[theorem]{Definition}

\newtheorem{lemma}[theorem]{Lemma}

\newtheorem{proposition}[theorem]{Proposition}
\newtheorem{remark}[theorem]{Remark}

\begin{document}

\date{}
\title{\textbf{
Local well-posedness for the great lake equation with transport noise}}
\author{Dan Crisan \qquad Oana Lang \\
\\ {\small Department of Mathematics, Imperial College London, SW7 2AZ, UK}}
\maketitle

\begin{abstract}
This work is a continuation of the authors' work in \cite{CL1}. In \cite{CL1} the equation satisfied by an incompressible fluid with stochastic transport is analysed. Here we lift the incompressibility constraint. Instead we assume a \textit{weighted} incompressibility condition. This condition is inspired by a physical model for a fluid in a basin with a free upper surface and a spatially varying bottom topography (see \cite{LevermoreOliverTiti2}). Moreover, we assume a different form of the vorticity to stream function operator that generalizes the standard Biot-Savart operator which appears in the Euler equation. These two properties are exhibited in the physical model called the \textit{great lake equation}. For this reason we refer to the model analysed here as the \textit{stochastic great lake equation}. Just as in \cite{CL1}, the deterministic model is perturbed with transport type noise. The new vorticity to stream function operator generalizes the $\boldsymbol{curl}$ operator and it is shown to have good regularity properties. We also show that the initial smoothness of the solution is preserved. The arguments are based on constructing a family of viscous solutions which is proved to be relatively compact and to converge to a truncated version of the original equation. Finally, we show that the truncation can be removed up to a positive stopping time.  \\

\vspace{7mm}
\vspace{7mm}

\vspace{1mm}
\end{abstract}

\section{Introduction} \label{introduction}

Consider the following two-dimensional stochastic system which models the evolution of the vorticity corresponding to an inviscid fluid on the two-dimensional torus $\mathbb{T}^2$: 
\begin{equation}\label{lake1}
d\omega_t + \mathcal{L}_{u_t} \omega_t dt + \displaystyle\sum_i \mathcal{L}_{\xi_i}\omega_t \circ dW_t^i = 0 
\end{equation} 
\begin{equation*}
\begin{aligned} 
&\nabla \cdot (bu) = 0 \\
&\omega = b^{-1} \boldsymbol{curl} \ v \\
& v = u + \frac{1}{6}\delta^2b^2\nabla(\nabla \cdot u).
\end{aligned} 
\end{equation*}
The initial condition $\omega _{0}$ is assumed to be an element of the weighted Sobolev space $\mathcal{W}_b^{k,2}(\mathbb{T}^2)$. The velocity of the fluid is denoted by $u_{t}$ and $\omega _{t}$ stands for the corresponding vorticity.  
The term $\mathcal{L}_{u_t}\omega_t$ can be interpreted as the Lie derivative which expresses the change of vorticity along the flow generated by the velocity vector field $u_t$, while $\mathcal{L}_{\xi_i}\omega_t$ is a circulation-preserving perturbation of this field (see \cite{CFH}, \cite{Holm2015}). A new variational approach for deriving stochastic partial differential equations which preserve fundamental properties of classical fluid dynamics has been introduced in \cite{Holm2015}. This method is now known as Stochastic Advection by Lie Transport (SALT). The model presented here is inspired by the so-called \textit{great lake equations} (see \cite{CamassaHolmLevermore1}, \cite{CamassaHolmLevermore2}, \cite{LevermoreOliverTiti2} for the deterministic case)  which model the circulation of an inviscid fluid in a shallow water basin with varying bottom topography $b$. The stochastic counterpart has been described in \cite{HolmLuesink} for a three-dimensional box. The domain we consider in this paper is the two-dimensional torus. The vector fields $\xi _{i}$ are time-independent and divergence-free and can be associated with uncertainty induced by missing physics or incomplete data (see \cite{Wei1}, \cite{Wei2}). The processes $W^{i}, i\geq 1$, are independent Brownian motions and $\delta$ is the aspect ratio of the domain (that is the ratio between vertical and horizontal length scales). The function $b$ is fixed in time and there exist two constants $b_{min}$ and $b_{max}$ such that $0 < b_{min} \leq b(x) \leq b_{max}$ for any $x \in \mathbb{T}^2$. If $b$ is constant then the great lake equations formally reduce to the classical 2D Euler equations (\cite{LevermoreOliverTiti2}). 

Nonlinear transport equations play a central role in modelling a broad range of phenomena such as flood waves, chemical reactions, gas dynamics, etc. Nonetheless, there are numerous small-scale physical processes which are still under-represented: turbulent multi-scale motion, friction, radiation, tropical cumulus convection (see e.g. \cite{Kalnay}). The introduction of a stochastic term of transport type aims at providing a better strategy for encoding randomness into a priori deterministic transport models. The noise structure is chosen such that the original physical properties of the model (e.g. Kelvin circulation theorem, energy conservation) are preserved also in the new stochastic setting: see \cite{Holm2015} or \cite{Memin}.

The stochastic part of our model follows the stochastic advection by Lie transport theory presented in \cite{Holm2015} and \cite{CS1}, and is represented by a stochastic integral of Stratonovich type. Well-posedness for the deterministic version of this system has been proven in \cite{LevermoreOliverTiti2}. In \cite{CL1} we have considered a similar stochastic vorticity equation but without taking into account the bottom topography $b$.  

 In this paper we prove that the system \eqref{lake1} admits a pathwise unique and probabilistically strong solution in the Sobolev space $\mathcal{W}_b^{k,2}(\mathbb{T}^2)$, the initial smoothness of the system being preserved. For a detailed analysis we use the It\^{o} form of the vorticity equation. Consequently, we must control a second order term which comes from the It\^{o} correction and involves the operator $\mathcal{L}_{\xi_i}^2$. We manage to do this by combining it with the quadratic variation of the stochastic integral. Further complications come out when trying to control the higher-order derivatives of the vorticity. Notwithstanding, we do this by proving a set of specific inequalities for the operator $\mathcal{L}_{\xi_i}$ (see Lemma \ref{apriori}) which are based on some smoothness and summability assumptions for the vector fields $\xi_i$ (see Assumptions \ref{xiassumpt}). These assumptions would be of particular interest when using this model as a signal process in stochastic filtering applications: see for instance \cite{Wei1}, \cite{Wei2}, \cite{Wei3} in the case of the stochastic 2D Euler equation. 

Note that the classical incompressibility condition $\nabla \cdot u = 0$ is not true here, a fact which usually generates technical difficulties when trying to control the nonlinear term and to obtain uniform a priori bounds. However, we have the weighted incompressibility condition $\nabla \cdot (bu) = 0$. Due to this, working in the weighted Sobolev space $\mathcal{W}_b^{k,2}(\mathbb{T}^2)$ comes as natural and we are able to extend the techniques from \cite{CL1} to this more general case. One of the key points here is to prove that the vorticity to stream function operator can be generalised (see Section \ref{biotsavartsection}) and therefore the smoothness of the velocity vector field is still controlled by the smoothness of its corresponding vorticity. We will show that even in a stochastic framework, $b$ affects just the geometry of the ambient space, not the topology, a positive consequence of the non-degeneracy of $b$. In \cite{LevermoreOliverTiti2} it has been shown that this property holds in a deterministic framework. In contrast to \cite{CL1} and following the idea from \cite{LevermoreOliverTiti2} we control the smoothness of the velocity vector field using Lax-Milgram-type arguments instead of Fourier analysis techniques.  

\vspace{3mm} 
In the sequel, $\mathbb{T}^2$ is the two-dimensional torus, $k\ge 0$ is a fixed positive integer and $\mathcal{W}_b^{k,2}$ is a weighted Sobolev space (see Section \ref{preliminaries}). 
The main result of this paper reads as follows: \\[2mm]
\textit{Theorem: Under certain conditions on the vector fields $(\xi_i)_i$, the vorticity equation of the two-dimensional stochastic great lake system \eqref{lake1} 
admits a pathwise unique local strong solution which belongs to the space $\mathcal{W}_b^{k,2}(\mathbb{T}^{2})$.} 
\begin{remark}\label{itoform}
The corresponding It\^{o} form of the above evolution equation for $\omega$ is
\begin{equation}\label{maineqnito}
d\omega_t + \mathcal{L}_{u_t} \omega_tdt + \displaystyle\sum_{i=1}^{\infty} \mathcal{L}_{\xi_i} \omega_tdW_t^i = \frac{1}{2}\displaystyle\sum_{i=1}^{\infty} \mathcal{L}_{\xi_i}^2 \omega_t dt.
\end{equation}
\end{remark}

The assumptions on the vector fields $(\xi_i)_i$ are described in Section \ref{preliminaries}. In brief, they are assumed to be sufficiently smooth and their corresponding norms to decay sufficiently fast as $i$ increases, so that the infinite sums in \eqref{lake1}, respectively, in \eqref{maineqnito}, make sense in the right spaces, (see condition \eqref{xiassumpt} below). \\
Since $u$ is such that $\nabla \cdot (bu) =0$, there exists (\cite{MajdaBertozzi}) a stream function $\psi$ such that $u = b^{-1}\nabla^{\perp}\psi$. According to \cite{CamassaHolmLevermore1}, using the weighted incompressibility condition $\nabla \cdot (bu) = 0$, the relationship between $v$ and $u$ is defined through a linear operator $\mathcal{M}$ as follows:

\begin{equation*}
\begin{aligned} 
v = \mathcal{M}u = u + \delta^2b^{-1}\bigg[-\frac{1}{3}\nabla(b^3\nabla \cdot u) - \frac{1}{2}\nabla(b^2 u \cdot \nabla b) + \frac{1}{2} b^2(\nabla \cdot u)\nabla b + b(u\cdot \nabla b)\nabla b \bigg].
\end{aligned} 
\end{equation*}
The introduction of a stream function here allows one to recover $u$ from $\omega$. Qualitatively, if $u = K\omega$ and $\omega$ is the solution of the vorticity equation, we want to show that $K$ has good regularity properties and that it imposes the right incompressibility conditions on $u$. Note that $\psi$ is the solution of the elliptic
problem
\begin{equation*}
\omega = b^{-1} \boldsymbol{curl} \big(\mathcal{M}(b^{-1}\nabla^{\perp}\psi)) \ \ \hbox{on $\mathbb{T}^2$}.
\end{equation*}
The great lake equation is more general than the Euler equation in the sense that the link between the vorticity field and the velocity field is not given by the $\boldsymbol{curl}$ operator only, but by a more general (linear) operator which is proved to have good regularity properties. 

The paper is organised as follows:  In Section \ref{preliminaries} we introduce the main assumptions and key notations. In Section \ref{mainresults} we introduce the main results. In Section \ref{uoriginal} we show that the solution of the great lake equation is almost surely pathwise unique. Section \ref{biotsavartsection} contains the regularity properties of the vorticity to stream function operator. In Section \ref{uetruncated} one can find the proof for existence of a strong solution (in the sense of Definition \ref{solutions}) and for uniqueness, of a truncated version of the great lake equation. In Section \ref{localtoglobal} we explain how the local solution of the great lake equation is obtained from the global solution of its truncated version. In Section \ref{contEuler} we show continuity with respect to initial conditions for the original equation. In Section \ref{sectiontightness} we show that the family of approximating solutions is relatively compact. In Section \ref{apriorisection} one can find a couple of technical properties which are essential in proving a priori estimates. \\

\section{Preliminaries}\label{preliminaries}

We summarise the notation used throughout the manuscript. Let $X$ be a generic Banach space and $b$ a \textit{weight function} as described above. Let $(\Omega ,\mathcal{F},(\mathcal{F}_{t})_{t\geq 0},\mathbb{P})$ be a filtered probability space with the sequence $(W^{i}) _{i\in\mathbb{N}}$ of
independent Brownian motions defined on it.
\begin{itemize}
\item We denote by $\mathbb{T}^{2}=\mathbb{R}^{2}/\mathbb{Z}^{2}$ the two-dimensional
torus.
\item $L_b^{p} (\mathbb{T}^{2}; X)$ is the class of all measurable $p$ - integrable functions $f$ defined on the two-dimensional torus, with values in $X$ ($p$ is a positive real number)\footnote{Here and later if the space $X$ coincides with the Euclidean space $\mathbb{R}$ or $\mathbb{R}^2$, it is omitted from the notation: For example $L^{p} (\mathbb{T}^{2}; X)$ becomes $L^{p} (\mathbb{T}^{2})$, etc.}. The space is endowed with its canonical weighted norm $\|f\|_{b,p} = \bigg(\displaystyle\int_{\mathbb{T}^2} \|f\|_X^p b(x)dx\bigg)^{1/p}$. Conventionally, for $p = \infty$ we denote by $L^{\infty}$ the space of essentially bounded measurable functions. 
\item For $a_1,a_2\in L_b^{2}\left(  \mathbb{T}^{2}\right)$, we denote by $\langle \cdot ,\cdot \rangle_b$ the
scalar product
\begin{equation*}
\langle a_1,a_2\rangle_b :=\displaystyle\int_{\mathbb{T}^{2}}a_1(x)\cdot a_2(x)b(x)dx.
\end{equation*}%
\item $\mathcal{W}_b^{m,p}(\mathbb{T}^{2})$ is the Sobolev space of functions $f \in L_b^p(\mathbb{T}^2)$ such that $D^{\alpha}f \in L_b^p(\mathbb{T}^2)$ for $0 \leq |\alpha| \leq m$, where $D^{\alpha}f$ is the distributional derivative of $f$ and $|\alpha|$ is the length of the multi-index $\alpha$. The canonical norm of this space is $\|f\|_{b,m,p} = \bigg(\displaystyle\sum_{0 \leq |\alpha| \leq m} \|D^{\alpha}f\|_{b,p}^p\bigg)^{1/p}$, with $m$ a positive integer and $1\leq p < \infty$. We denote by $\mathcal{W}_b^{-m,p}$ the dual of $\mathcal{W}_b^{m,p}$. When $p=\infty$ we define $\|f\|_{k,\infty}:= \displaystyle\max_{|\alpha| \leq k}\|D^{\alpha}f\|_{\infty}=\displaystyle\max_{|\alpha|\leq k}(ess \displaystyle\sup_x|D^{\alpha}f(x)|)$. A detailed presentation of Sobolev and weighted Sobolev spaces can be found in \cite{Adams} and \cite{KufnerWeightedSobSpaces}, respectively. 

\item $C^m(\mathbb{T}^2; X)$ is the (vector) space of all $X$-valued functions $f$ which are continuous on $\mathbb{T}^2$ with continuous partial derivatives $D^{\alpha}f$ of orders $|\alpha| \leq m$, for $m \geq 0$. $C^{\infty}(\mathbb{T}^2; X)$ is regarded as the intersection of all spaces $C^{m}(\mathbb{T}^2; X)$. 
\item $C([0, \infty); X)$ is the space of continuous functions from $[0, \infty)$ to $X$ equipped with 
the uniform convergence topology over compact subintervals of $[0, \infty)$.

\item $D([0, \infty); X)$ is the space of c\`{a}dl\`{a}g functions that is functions $f:[0, \infty) \rightarrow X$ which are right-continuous and have limits to the left, endowed with the Skorokhod topology. This topology is a natural choice in this case because its corresponding metric transforms $D([0, \infty); X)$ into a complete separable metric space. For further details see \cite{EthierKurtz} Chapter 3, Section 5, pp. 117-118. 
\item Given $a_1: \mathbb{T}^{2}\rightarrow\mathbb {R}^2$, we define the differential
operator $\mathcal{L}_{a_1}$ by 
$\mathcal{L}_{a_1}a_2:=a_1\cdot \nabla a_2$
 for any map $a_2: \mathbb{T}^{2}\rightarrow\mathbb {R}$ such that the weighted scalar product
 between $a_1$ and $a_2$ makes sense. In line with this, 
\begin{equation*} 
\mathcal{L}_{i}\omega _{t}:=\mathcal{L}_{\xi _{i}}\omega _{t}:=\xi
_{i}\cdot \nabla \omega _{t}\ \ \ \hbox{and}\ \ \ \mathcal{L}%
_{i}^{2}\omega _{t}:=\mathcal{L}_{\xi _{i}}^{2}\omega _{t}:=\xi _{i}\cdot
\nabla (\xi _{i}\cdot \nabla \omega _{t}).
\end{equation*}%
Denote the dual of $\mathcal{L}$ by $\mathcal{L}^{\star}$ that is $%
\langle \mathcal{L} a_1,a_2\rangle =\langle a_1,\mathcal{L}^{\star }a_2\rangle$. 
\item $\mathbb{D}_1:= \left\{ u \in C^{\infty}(\mathbb{T}^2): \nabla \cdot (bu) = 0\right\}$.
\item $\mathbb{D}_2 := $ the completion of $\mathbb{D}_1$ in the $L_b^2$ - norm. 
\item $\mathbb{D}_3 := $ the completion of $\mathbb{D}_1$ in the $\mathcal{W}_b^{1,2}$ - norm.
\item $\mathcal{P}:=$ the $L_b^2$ - orthogonal projection onto $\mathbb{D}_2$.  
\item For any vector $u \in \mathbb{R}^2$ we denote the gradient of $u$ by $\nabla u = (\partial_1u, \partial_2 u)$ and the corresponding orthogonal by $\nabla^{\perp}u = (\partial_2 u, -\partial_1 u )$.
\end{itemize}

\begin{remark}\label{dual}
If $\boldsymbol{div} \left(b\xi _{i}\right)=\nabla \cdot \left(b\xi _{i}\right)=0$, then the dual of the operator $\mathcal{L}_{i} $ is $-\mathcal{L}_{i}$.
\end{remark}
\noindent\textbf{Assumptions on the vector fields $(\xi_i)_i$} \newline
The vector fields $\xi_i:\mathbb{T}^2 \rightarrow \mathbb{R}%
^2 $ are chosen to be time-independent
quantities such that $\nabla \cdot (b\xi_i)=0$, which need to be specified from the underlying physics. We assume that for any $f \in \mathcal{W}_b^{2,2}(\mathbb{T}^2)$:
\begin{equation} \label{xiassumpt}
\begin{split} 
 \displaystyle\sum_{i=1}^{\infty} \|\mathcal{L}_{i}f\|_{b,2}^2 \leq C\|f\|_{b,1,2}^2, \ \ \ \ \ \ \ 
\displaystyle\sum_{i=1}^{\infty} \|\mathcal{L}_i^2f\|_{b,2}^2 \leq C\|f\|_{b,2,2}^2, \ \ \ \ \ \ \
\displaystyle\sum_{i=1}^{\infty} \|\xi_i\|_{k+1,\infty}^2 < \infty.
\end{split}
\end{equation}
Provided that $\omega\in L_b^2(0, T; \mathcal{W}_b^{2,2}(\mathbb{T}^2, \mathbb{R}))$, the first two conditions in \eqref{xiassumpt} ensure that the infinite sums of stochastic integrals
\begin{equation} \label{infsum1}
\sum_{i=1}^{\infty} \int_0^t \xi_i \cdot \nabla\omega_sdW_s^i 
\end{equation}
are well defined and belong to $L_b^2(0, T; L_b^2(\mathbb{T}^2, \mathbb{R}))$
Similarly,  the processes $s\rightarrow\xi_i \cdot \nabla(\xi_i \cdot \nabla\omega_s)$
are well-defined and belong (pathwise) to $L_b^2(0, T; L_b^2(\mathbb{T}^2, \mathbb{R}))$. In particular, the It\^o correction in \eqref{maineqnito} is well-defined.
The third condition is needed for proving a number of required a priori estimates (see Lemma \ref{aprioriapprox}).
In the following definitions $k\ge 2$ is a fixed integer.  
\begin{definition}\label{solutions}$\left.\right.$\\[-5mm]
\begin{enumerate}
\item[a.] A strong\footnote{Note that this solution is strong in both probabilistic and PDE sense, as it is defined on the probability space fixed in advance.} solution of the
stochastic partial differential equation \eqref{maineqnito} is an $(\mathcal{F}_{t})_{t}$-adapted process $\omega :\Omega \times \mathbb{T}^{2}\rightarrow \mathbb{R}$ with trajectories in the space $C([0,\infty);\mathcal{W}_b^{k,2}(\mathbb{T}^{2}))$, such that the identity 
\begin{equation*}
\omega _{t}=\omega _{0}-\int_{0}^{t}u_{s}\cdot \nabla \omega
_{s}ds-\sum_{i=1}^{\infty }\int_{0}^{t}\xi _{i}\cdot \nabla \omega
_{s}dW_{s}^{i}+\frac{1}{2}\sum_{i=1}^{\infty }\int_{0}^{t}\xi _{i}\cdot
\nabla \big(\xi _{i}\cdot \nabla \omega _{s}\big)ds
\end{equation*}%
with $\omega _{|_{t=0}}=\omega _{0}$, holds $\mathbb{P}$-almost surely  in $L_b^{2}(\mathbb{T}^{2};\mathbb{R})$.
\item[b.] A weak/distributional solution of
equation \eqref{maineqnito} is an $(\mathcal{F}_t)_t$-adapted process $\omega : \Omega
\times \mathbb{T}^2 \rightarrow \mathbb{R}$ with trajectories in the set $%
C([0, \infty); L_b^{2}(\mathbb{T}^2))$, which satisfies the
equation  i.e. 
\begin{equation*}
\langle\omega_t, \varphi\rangle_b = \langle\omega_0,\varphi\rangle_b- %
\displaystyle\int_0^t \langle \omega_s, \mathcal{L}_{u_s}^{\star}\varphi
\rangle_b ds - \displaystyle\sum_{i=1}^{\infty} \displaystyle\int_0^t
\langle\omega_s, \mathcal{L}_{i}^{\star}\varphi\rangle_b dW_s^i + \frac{1}{2}%
\displaystyle\sum_{i=1}^{\infty}\displaystyle\int_0^t\langle\omega_s, %
\mathcal{L}_{i}^{\star}\mathcal{L}_{i}^{\star}\varphi\rangle_b ds
\end{equation*}
holds $\mathbb{P}$-almost surely for all $\varphi\in C^{\infty}(\mathbb{T}%
^2, \mathbb{R})$.
\item[c.] A weak probabilistic solution of equation \eqref{maineqnito} is a triple $(\tilde\omega, (\tilde W^i)_i), (\tilde\Omega, \mathcal{\tilde F}, \tilde{\mathbb{P}}), (\mathcal{\tilde F}_t)_t$ such that $(\tilde\Omega, \mathcal{\tilde F}, \tilde{\mathbb{P}})$ is a probability space, $(\mathcal{\tilde F}_t)_t$ is a filtration defined on this space,   
$\tilde\omega : \Omega
\times \mathbb{T}^2 \rightarrow \mathbb{R}$ is a continuous $(\mathcal{\tilde F}_t)_t$-adapted process with trajectories in the set $%
C([0, \infty); \mathcal{W}_b^{k,2}(\mathbb{T}^2))$, $(\tilde W^i)_i$ are independent $(\mathcal{\tilde F}_t)_t$-adapted Brownian motions, and the identity 
\begin{equation*}
\omega _{t}=\omega _{0}-\int_{0}^{t}u_{s}\cdot \nabla \omega
_{s}ds-\sum_{i=1}^{\infty }\int_{0}^{t}\xi _{i}\cdot \nabla \omega
_{s}dW_{s}^{i}+\frac{1}{2}\sum_{i=1}^{\infty }\int_{0}^{t}\xi _{i}\cdot
\nabla \big(\xi _{i}\cdot \nabla \omega _{s}\big)ds
\end{equation*}%
with $\omega _{|_{t=0}}=\omega _{0}$, holds $\mathbb{\tilde{P}}$-almost surely in  $L_b^{2}(\mathbb{T}^{2};\mathbb{R})$. 
\item[d.]A classical solution of equation \eqref{maineqnito}
is an $(\mathcal{F}_t)_t$-adapted process $\omega : \Omega \times [0, \infty) \times \mathbb{T%
}^2 \rightarrow \mathbb{R}$ with trajectories of class $C([0, \infty); C^2(%
\mathbb{T}^2;\mathbb{R}))$.

\item[e. ] A local solution of equation \eqref{maineqnito} is given by a pair $(\omega, \tau)$ consisting of a stopping time $\tau : \Omega \rightarrow [0,\infty)$ and a process $\omega: \Omega \times [0,\tau]\times\mathbb{T}^2\rightarrow \mathbb{R}$ such that the solution trajectory is of class $C([0,\tau];\mathcal{W}_b^{k,2}(\mathbb{T}^2))$, $\omega_{t \wedge \tau}$ is $(\mathcal{F}_t)_t$-adapted, and for any stopping time $\tau'\leq \tau$ the integral identity
\begin{equation*}
    \omega _{\tau'}=\omega _{0}-\int_{0}^{\tau'}u_{s}\cdot \nabla \omega
_{s}ds-\sum_{i=1}^{\infty }\int_{0}^{\tau'}\xi _{i}\cdot \nabla \omega
_{s}dW_{s}^{i}+\frac{1}{2}\sum_{i=1}^{\infty }\int_{0}^{\tau'}\xi _{i}\cdot
\nabla \big(\xi _{i}\cdot \nabla \omega _{s}\big)ds
\end{equation*}
with $\omega _{|_{\tau'=0}}=\omega _{0}$, holds $\mathbb{P}$-almost surely in  $L_b^{2}(\mathbb{T}^{2};\mathbb{R})$.

\end{enumerate}
\end{definition}
\begin{remark}\label{whatisu} The velocity field $v$ is not uniquely identified through the equation $\omega=b^{-1}\boldsymbol{curl}\  v$. Indeed any two velocity fields that differ by a constant will lead to the same vorticity map $\omega$. Instead we identify $v$ through the "explicit" formula $v=\mathcal{M}\left( b^{-1}\nabla^{\perp}\psi \right)$, where $\psi$ is a stream function which exists due to the fact that $\nabla \cdot (bu) = 0$.         
\end{remark} 

\begin{remark}
Note that $\omega_t \in \mathcal{W}_b^{k,2}(\mathbb{T}^2)$ implies $u_t \in \mathcal{W}_b^{k+1,2}(\mathbb{T}^2)$ (see Section \ref{biotsavartsection}). By standard Sobolev embedding theorems $\mathcal{W}_b^{k+1,2}(\mathbb{T}^2) \hookrightarrow \mathcal{W}_b^{k,2}(\mathbb{T}^2) \hookrightarrow C(\mathbb{T}^2)$ for $k\geq 2$, hence the terms $\mathcal{L}_{u_t}\omega_t = u_t \cdot \nabla \omega_t\in L_b^2(\mathbb{T}^2, \mathbb{R})$ and $\langle \omega_s, \mathcal{L}_{u_s}^{\star}\varphi
\rangle_b$ are well defined. 
\end{remark}
\begin{remark} \label{weak} Naturally, if $\omega_t$ is a strong solution in
the sense of Definition \ref{solutions}a, then it is also a\\ weak/distributional solution in the sense of Definition \ref{solutions}b. 
Note also that if $\omega_t$ is a weak/distributional solution with paths in $C([0,T]; \mathcal{W}_b^{k,2}(\mathbb{T}^2))$ then, by integration by parts, it is also a (probabilistic) strong solution. 
\end{remark}

\section{Main results}\label{mainresults}
Let $\left (\Omega, \mathcal{F}, (\mathcal{F}_t)_t, \mathbb{P} \right)$ be a filtered probability space and $k \geq 2$ a fixed integer as above. The main result of the paper reads as follows: 
\begin{theorem} \label{mainresult1}
If $\omega_0 \in \mathcal{W}_b^{k,2}(\mathbb{T}^2)$ then the two-dimensional stochastic great lake equation \eqref{lake1}
\begin{equation*}
d{\omega}_t+ \mathcal{L}_{u_t}{\omega}_t dt + \displaystyle %
\sum_{i=1}^{\infty}\mathcal{L}_{i} ^{} \omega_t\circ dW_t^i=0  \label{1}
\end{equation*}
admits a pathwise unique and probabilistically strong local solution (in the sense of Definition \ref{solutions})
$\omega=\{\omega_t,t\in [0,\infty)\}$ with trajectories in the space 
$C( [  0,\infty);\mathcal{W}_b^{k,2}(\mathbb{T}^2))$. In particular,
if  $k\geq 4$ the solution is classical. 
\end{theorem}

Moreover, we have continuity with respect to initial conditions:    
\begin{theorem}\label{mainresult2}
Let $\omega$, $\tilde\omega$ be two $\mathcal{W}_b^{k,2}(\mathbb{T}^2)$-solutions of the great lake equation \eqref{lake1}. Define the process $B = (B_t)_t$ as $B_t:=\displaystyle\int_0^t\|\omega_s\|_{b,k,2}ds$, for any $ t\ge 0$. Then there exists a positive constant $\mathcal{C}$ independent of the two solutions $\omega$ and $\tilde{\omega}$ such that 
  
\begin{equation}\label{conteqn}
\mathbb E [e^{-\mathcal{C}B_t}||\omega_t-\tilde\omega_t||_{b,k-1,2}^2]\le 
||\omega_0-\tilde\omega_0||_{b,k-1,2}^2.
\end{equation}
\end{theorem}
The proof of these two theorems are covered in Sections \ref{uoriginal}, \ref{localtoglobal}, and \ref{uetruncated}.

\section{Pathwise uniqueness} \label{uoriginal}
In this section we prove that any two solutions $\omega^1$ and $\omega^2$ defined on the same probability space $\left( \Omega, \mathcal{F}, \mathbb{P} \right)$, driven by the same Brownian motion $(W^i)^i$, and with $\mathbb{P}$-almost surely the same initial conditions $\omega_0^1$ and $\omega_0^2$, are indistinguishable, that is 
\begin{equation*}
\mathbb{P}\left(\omega_t^1 = \omega_t^2  \ \hbox{for all} \  t \geq 0\right) = 1. 
\end{equation*} 
From a probabilistic perspective this means that the solution is pathwise unique. This result is essential for the proof of existence of a strong solution in the sense of Definition \ref{solutions}.  
\\

\noindent Suppose that equation \eqref{lake1} admits two $\mathcal F_t$-adapted solutions $\omega^1$ and $%
\omega^2 $ with trajectories in the space 
$C([  0,\infty);\mathcal{W}_b^{k,2}(\mathbb{T}^2))  $ and let $\bar{\omega}:=\omega^1-\omega^2$. Consider the
corresponding velocities $u^1$ and $u^2$ such that $b^{-1}\boldsymbol{curl} \left(\mathcal{M} u^1\right)=\omega^1$, $b^{-1}\boldsymbol{curl} \left(\mathcal{M} u^2\right)=\omega^2$ and $\bar{u}:=u^1-u^2$. Since both $\omega^1$ and $%
\omega^2$ satisfy \eqref{maineqnito}, their difference satisfies 
\begin{equation*}
d{\bar{\omega}}_t+(\mathcal{L}_{\bar{u}_t}\omega_t^1 + \mathcal{L}_{u_t^2} {\bar{\omega}_t}) dt + \displaystyle \sum_{i=1}^{\infty}\mathcal{L}_{i}{%
\bar{\omega}_t} dW_t^i-\frac{1}{2}\displaystyle\sum_{i=1}^{\infty}\mathcal{L}_{i}^2 \bar{\omega}_tdt=0.
\end{equation*}
By the It\^{o} formula one obtains 
\begin{equation*}
\begin{aligned}
d\|\bar{\omega}_t\|_{b,2}^2&=-2\displaystyle\sum_{i=1}^{\infty}\langle
\bar{\omega}_t, \mathcal{L}_{i}\bar{\omega}_t \rangle_b dW_t^i-2\langle
\bar{\omega}_t, \mathcal{L}_{\bar{u}_t}\omega_t^1 + \mathcal{L}_{u_t^2} {\bar{\omega}_t}\rangle_b dt \\ &+ \displaystyle\sum_{i=1}^{\infty}
\left(\big\langle \bar{\omega}_t, \mathcal{L}_{i}^2
\bar{\omega}_t \big\rangle_b  + \langle\mathcal{L}_{i} \bar{\omega}_t, \mathcal{L}_{i} \bar{\omega}_t \rangle_b\right)
dt. \end{aligned}
\end{equation*}
Note that the last term in the above identity is null (see Lemma \ref{apriori}) 
and that
\vspace{-3mm} 
\begin{equation*}
 |\langle
\bar{\omega}_t, \mathcal{L}_{\bar{u}_t}\omega_t^1 \rangle_b| \leq
\|\bar{\omega}_t\|_{b,2}\|\bar{u}_t\|_{b,4}\|\nabla\omega_t^1\|_{b,4} \leq
C\|\bar{\omega}_t\|_{b,2}^2\|\omega_t^1\|_{b,k,2}.
\end{equation*}
This is true since by the Sobolev embedding theorem (see \cite{Adams} Theorem 4.12 case A) 
one has $\|\nabla\omega_{t}^1\|_{b,4} \leq C \|\omega_t^1\|_{b,k,2}$ and due to the smoothness properties proved in Section \ref{biotsavartsection} one has
$\|\bar{u}_t\|_{b,4} \leq C\|\bar{u}_t\|_{b,1,2} \leq C \|\bar{\omega}_t\|_{b,2}$. 

Finally, observe that 
$\langle
\bar{\omega}_t,  \mathcal{L}_{u_t^2} {\bar{\omega}_t}\rangle_b =-\frac{1}{2} \displaystyle\int_{\mathbb{T}^2} (\nabla \cdot (bu_t^2))
(\bar{\omega}_t)^2 dx =0 $
since $\nabla \cdot(bu_t^2)=0$. It follows that 
\begin{equation}\label{lastinequality}
\begin{aligned} d\|\bar{\omega}_t\|_{b,2}^2&=-2\langle \bar{\omega}_t, \mathcal{L}_{\bar{u}_t}\omega_t^1 \rangle_b dt \leq
C\|\omega_t^1\|_{b,k,2}\|\bar{\omega}_t\|_{b,2}^2dt. \end{aligned}
\end{equation}
Since we only have a priori bounds for the expected value of the process $t\rightarrow \|\omega_t^1\|_{b,k,2}$ and not for its pathwise values, the uniqueness cannot be deduced through a classical Gronwall-type argument. Instead, we proceed as follows: let $B = (B_t)_t$ be the process defined as 
$B_t:=\displaystyle\int_0^tC\|\omega_s^1\|_{b,k,2}ds$, for any $ t\ge 0$.
This is an increasing process that stays finite $\mathbb{P}$-almost surely for all $t\ge 0$ as the paths of $\omega^1$ are in $C( [  0,\infty);\mathcal{W}_b^{k,2}(\mathbb{T}^2))  $. By the product rule, 
\begin{equation*}
d\big(e^{-B_t}\|\bar{\omega}_t\|_{b,2}^2\big)=e^{-B_t}(d\|\bar{\omega}_t\|_{b,2}^2-C\|\bar{\omega}_t\|_{b,2}^2\|%
\omega_t^1\|_{k,2}dt)\leq 0,
\end{equation*}
which leads to \vspace{-3mm} 
\begin{equation*}
\begin{aligned} e^{-B_t}\|\bar{\omega}_t\|_{b,2}^2 &\leq 0.
\end{aligned}
\end{equation*}
We conclude that 
$e^{-B_t}\|\bar{\omega}_t\|_{2}^2=0$, and since $e^{-B_t}$ cannot be null due to the finiteness of $B_t$ we deduce that $\|\bar{\omega}_t\|_{b,2}^2=0$ almost surely, which gives the claim.

\section{The vorticity to stream function map}\label{biotsavartsection}

In this section we provide a generalisation of the Biot-Savart law from \cite{CL1}, which is the main ingredient used to extend the well-posedness properties to a more general class of stochastic partial differential equations. Note that the stochastic part does not interfere with this generalisation. This is due to the fact that stochasticity has been introduced as a constraint within a variational principle such that the physical meaning of the quantities of interest has not been altered: see \cite{Holm2015}, \cite{CS1}. More precisely, $\omega = b^{-1} \boldsymbol{curl} \ v$ where $v$ depends on the deterministic, unperturbed velocity vector field $u$ only, without any stochasticity involved. The central result is Proposition \ref{biotsavart} \footnote{The proofs from this section follow closely the arguments from \cite{LevermoreOliverTiti2} Section 2.}. 

\begin{lemma}
The restriction of the operator $\mathcal{M}$ to the space $\mathbb{D}_2$ is a continuous, positive, invertible and self-adjoint operator such that for every $u \in L_b^2(\mathbb{T}^2)$ one has
\begin{equation*}
\mathcal{P}\mathcal{M}\mathcal{P}u = \mathcal{P}\left(u + \frac{\delta^2}{3}\nabla b \left(\nabla b \cdot \mathcal{P}u\right)\right). 
\end{equation*} 
\end{lemma}

\noindent \textbf{Proof:}
Let $u,v \in \mathbb{D}_1$. We have
\begin{equation*}
\begin{aligned}
\mathcal{M}u = u + \delta^2 \left( - \frac{1}{3}b^{-1} \nabla\left(b^3 \nabla \cdot u\right) - \frac{1}{2}b^{-1} \nabla \left( b^2 u \cdot \nabla u\right) + \frac{1}{2} b \left(\nabla \cdot u \right)\nabla b + \left(u \cdot \nabla b \right)\nabla b\right)
\end{aligned}
\end{equation*}
and the bilinear form $\ell$ considered below can be written as 
\begin{equation*}
\begin{aligned}
\ell(u,v) & := \langle \mathcal{M}u, v\rangle_b = \displaystyle\int_{\mathbb{T}^2} \mathcal{M}u \cdot v bdx \\
& = \langle u,v\rangle_b  - \frac{\delta^2}{3}\displaystyle\int_{\mathbb{T}^2}b^{-1} \nabla\left(b^3 \nabla \cdot u\right) \cdot vb dx - \frac{\delta^2}{2}\displaystyle\int_{\mathbb{T}^2}b^{-1} \nabla \left( b^2 u \cdot \nabla u\right) \cdot v bdx \\
& + \frac{\delta^2}{2} \displaystyle\int_{\mathbb{T}^2} b^2(\nabla \cdot u)(\nabla b \cdot v)dx + \delta^2 \displaystyle\int_{\mathbb{T}^2} b(u \cdot \nabla b)(\nabla b \cdot v)dx \\
& = \langle u,v\rangle_b + \frac{\delta^2}{3}\displaystyle\int_{\mathbb{T}^2} b^3(\nabla \cdot u)(\nabla \cdot v)dx + \frac{\delta^2}{2} \displaystyle\int_{\mathbb{T}^2} b^2(u \cdot \nabla b)(\nabla \cdot v)dx \\
& + \frac{\delta^2}{2}\displaystyle\int_{\mathbb{T}^2} b^2 (\nabla \cdot u)(v \cdot \nabla v)dx + \delta^2\displaystyle\int_{\mathbb{T}^2} b(u \cdot \nabla b)(v \cdot \nabla b)dx \\
& = \langle u,v \rangle_b - \frac{\delta^2}{3} \displaystyle\int_{\mathbb{T}^2} b (u \cdot \nabla b)(v \cdot \nabla b) dx - \frac{\delta^2}{2}\displaystyle\int_{\mathbb{T}^2} b (u \cdot \nabla b)(v \cdot \nabla b)dx \\
& - \frac{\delta^2}{2}\displaystyle\int_{\mathbb{T}^2} b(u \cdot \nabla b)(v \cdot \nabla b)dx + \delta^2 \displaystyle\int_{\mathbb{T}^2} b (u \cdot \nabla b)(v \cdot \nabla b)dx \\
& = \langle u,v \rangle_b + \frac{\delta^2}{3}\langle u \cdot \nabla b, v \cdot \nabla b\rangle_b \\
& \leq c\|u\|_{b,2}\|v\|_{b,2}. 
\end{aligned}
\end{equation*}
The calculations have been simplified by the fact that $u,v \in \mathbb{D}_1$ implies $\nabla \cdot u = -b^{-1}u \cdot \nabla b$ and $\nabla \cdot v = -b^{-1}v \cdot \nabla b$. Therefore $\ell$ is continuous and symmetric and due to the Poincar\'{e} inequality it is also coercive. By the Lax-Milgram theorem (see for instance \cite{Brezis}, pp. 140, Corollary 5.8) we can conclude that the operator $\mathcal{M}$ is invertible. 

\begin{proposition}\label{biotsavart}
For every $\omega \in \mathcal{W}_b^{-1,2}(\mathbb{T}^2)$ there exists a unique function $u \in \mathbb{D}_2$ such that $u = K\omega$. Moreover, $K$ is continuous and 
\begin{equation*}
\|K\omega\|_{b,k,p} \leq C\|\omega\|_{b,k-1,p}.
\end{equation*}
\end{proposition}  

\noindent\textbf{Proof}: Let $u, v \in \mathbb{D}_1$. Then (see \cite{MarchioroPulvirenti}) there exist $\varphi_1, \varphi_2 \in C^{\infty}(\mathbb{T}^2)$ such that $bu = \nabla^{\perp} \varphi_1$ and $bv = \nabla^{\perp} \varphi_2$. We integrate by parts and use the previous lemma to write
\begin{equation*}
\begin{aligned} 
\tilde{\ell}(\varphi_1, \varphi_2)&:= \langle \varphi_1, \omega \rangle_b = 
\left\langle \varphi_1, b^{-1} \boldsymbol{curl}\left(\mathcal{M}\left(b^{-1}\nabla^{\perp}\varphi_2\right)\right) \right\rangle_b \\
 & = \displaystyle\int_{\mathbb{T}^2} \left( b^{-1}\nabla^{\perp}\varphi_1\right) \cdot \mathcal{M} \left( b^{-1} \nabla^{\perp}\varphi_2 \right) bdx \\
 & \leq C \|u\|_{b,2}\|v\|_{b,2}. 
\end{aligned}
\end{equation*}
On the other hand, we can also write
\begin{equation*}
\left\langle \varphi_1, b^{-1} \boldsymbol{curl}\left(\mathcal{M}\left(b^{-1}\nabla^{\perp}\varphi_2\right)\right)\right\rangle_b \leq C\|\varphi_1\|_{b,1,2}\|\varphi_2\|_{b,1,2}
\end{equation*}
and by the Poincar\'{e} inequality
\begin{equation*}
\|\varphi_1\|_{b,1,2}^2 \leq C \left\langle \varphi_1, b^{-1} \boldsymbol{curl}\left(\mathcal{M}\left(b^{-1}\nabla^{\perp}\varphi_1\right)\right)\right\rangle_b
\end{equation*}
which explains also the elliptic character of the operator. 
The same estimates hold for $\varphi_1, \varphi_2 \in \mathcal{W}^{1,2}(\mathbb{T}^2)$ by a density argument. 
Therefore the bilinear form $\tilde{\ell}$ is continuous and coercive. By the Lax-Milgram theorem, for any $\omega \in \mathcal{W}^{-1,2}(\mathbb{T}^2)$ there exists a unique $u \in \mathbb{D}_2$ such that $u = K\omega$ and $K: \mathcal{W}^{-1,2}(\mathbb{T}^2) \rightarrow \mathbb{D}_2$ is continuous. Note that the problem
\begin{equation*}
\omega = b^{-1} \boldsymbol{curl} \big(\mathcal{M}(b^{-1}\nabla^{\perp}\psi)) \ \ \hbox{on $\mathbb{T}^2$}
\end{equation*}
is elliptic and therefore by regularity theorems for elliptic problems (see \cite{Evans} Section 6.3 or \cite{Warner} Section 6.2.8) we actually have
\begin{equation*}
\|\psi\|_{b,k,p} \leq C \|\omega\|_{b,k-2,p}. 
\end{equation*}

\section{Local existence for the great lake equation}\label{localtoglobal}
The existence of the solution of equation \eqref{maineqnito} is proved by first showing that a truncated version of it has a solution, and then removing the truncation up to a positive stopping time. In particular we truncate the nonlinear term in \eqref{maineqnito} by using a smooth function $f_R$ equal to $1$ on $[0,R]$, equal to $0$ on $[R+1, \infty)$, and decreasing on $[R, R+1]$, for arbitrary $R>0$, with $f_R(u_t^R) := f_R(\|u_t^R\|_{b,k,2})$. 
Then we have the following: 
\begin{proposition} \label{truncatedEulert}
If $\omega_0 \in \mathcal{W}_b^{k,2}(\mathbb{T}^2)$ such that $\nabla \cdot (b\omega_0) = 0$, then the following equation 
\begin{equation}\label{truncatedEulere}
d{\omega}_t^{R}+ f_R(u_t^R)\mathcal{L}_{u_t^R}{\omega}_t^R dt + \displaystyle %
\sum_{i=1}^{\infty}\mathcal{L}_{i} ^{} \omega_t^R\circ dW_t^i=0  
\end{equation}
admits a  unique global $\mathcal F_t$-adapted solution
$\omega^R=\{\omega^R_t,t\in [0,\infty)\}$ with trajectories in the space 
$C( [  0,\infty);\mathcal{W}_b^{k,2}(\mathbb{T}^2))  $. In particular,
if  $k\geq 4$, the solution is classical. 
\end{proposition}

\begin{remark}
Observe that, by definition, the truncation function $f_R$ depends on the norm $\|\omega_t^{R}\|_{b,k-1,2}$ and not on the norm $\|\omega_t^{R}\|_{b,k,2}$. This is not incidental as it suffices to control the norm $\|u_t^{R}\|_{b,k,2}$  (see Section \ref{biotsavartsection}).    
\end{remark}
We prove Proposition \ref{truncatedEulert} in Section \ref{uetruncated}. For now let us proceed with the proof of local existence for the solution of the stochastic great lake equation \eqref{maineqnito}. Define the stopping time 
$$\tau_R(\omega):=\inf_{t\ge 0}\left\{\|\omega_t^{R}\|_{b,k-1,2}\geq \frac{R}{\mathcal{C}} \right\}$$
where $\mathcal{C}$ is such that $\|\nabla u\|_{\infty} \leq \mathcal{C} \|\omega\|_{b,k,2}$. Observe that such a constant exists due to the Sobolev embedding $\mathcal{W}_b^{k,2} \hookrightarrow L^{\infty}$ and to the regularity properties proven in Section \ref{biotsavartsection}, since
$\|\nabla u\|_{\infty} \leq \mathcal{C}\|\nabla u\|_{b,k,2} \leq \mathcal{C}\|u\|_{b,k+1,2} \leq \mathcal{C}\|\omega\|_{b,k,2}$.
\begin{lemma}\label{localslnjustif}
Let $\omega_0 \in \mathcal{W}_b^{k,2}(\mathbb{T}^2)$ and $\omega^R : \Omega \times [0, \infty) \times \mathbb{T}^2 \rightarrow \mathbb{R}$ be a global $\mathcal{W}_b^{k,2}(\mathbb{T}^2)$-solution of the truncated equation \eqref{truncatedEulere} and $\omega:\Omega \times [0, \tau_R] \times \mathbb{T}^2 \rightarrow \mathbb{R}$ be the restriction of $\omega^R$ to the time interval $[0, \tau_R]$. Then $\omega$ is a local $\mathcal{W}_b^{k,2}(\mathbb{T}^2)$-solution of the original great lake equation \eqref{maineqnito}.
\end{lemma}
\noindent\textbf{\textit{Proof}}
Observe that for $t \in [0, \tau_R]$ we have $\|\nabla u_t\|_{\infty} \leq \mathcal{C} \|\omega_t\|_{b,k,2} \leq R$. Thus, $f_R(u_t^R)=1$ and taking into account the pathwise uniqueness property we conclude that the truncated equation coincides with the original equation. 
\begin{remark}
If we dispense with the requirement of showing the existence of a strong solution, then similar uniqueness and relative compactness arguments can be used to show the existence of a unique \underline{global} weak solution, provided $\omega_0 \in L^{\infty}(\mathbb{T}^2)$. The key result here is the fact that the $L^{\infty}$-norm of the solution of \eqref{maineqnito} as well as that of any of its truncated versions remains constant in time (see Lemma \ref{apriori}):
\begin{equation*}
\begin{aligned} 
&\|\omega_t\|_{\infty} = \|\omega_0\|_{\infty} \\
&\|\omega_t^R\|_{\infty} \leq \|\omega_0^R\|_{\infty}.
\end{aligned} 
\end{equation*}
This is an extension of the result in \cite{Brzezniak}. In particular, one can show the existence of a global solution $\omega_t \in \mathcal{W}_b^{-2,2}(\mathbb{T}^2)$. 
\end{remark}

\section{Uniqueness and existence for the truncated equation}\label{uetruncated}
In this section we prove that the truncated equation admits a global unique solution. Any global solution for the truncated equation, restricted to the corresponding stopping time, is a local solution for the original great lake equation. 
\subsection{Pathwise uniqueness for the truncated equation} \label{utruncated}
The strategy for showing that any solution of the truncated equation is pathwise unique is similar to the one presented in Section \ref{uoriginal}. The only difference arises due to the truncated terms. We have 
\begin{equation*}
d\|\bar{\omega}_t\|_{b,2}^2+2\displaystyle\sum_{i=1}^{\infty}\langle
\bar{\omega}_t, \mathcal{L}_{i}\bar{\omega}_t \rangle_b dW_t^i=-2\langle
\bar{\omega}_t,(f_R(\omega_t^1)\mathcal{L}_{u^1_t}-f_R(\omega_t^2)\mathcal{L}_{u^2_t})\omega_t^1  \rangle_b dt  
 \end{equation*}
One can show that (see \cite{CFH} for a proof) there exists a constant $C=C(R)$ such that  
\[
\|f_R(\omega_t^1)u^1_t-f_R(\omega_t^2){u^2_t}\|_{b,4}\le C\|\bar \omega_t\|_{b,k-1,2}
\]
and therefore 
\begin{equation*}
 |\langle
\bar{\omega},(f_R(\omega_t^1)\mathcal{L}_{u^1_t}-f_R(\omega_t^2)\mathcal{L}_{u^2_t})\omega_t^1  \rangle_b| \leq
C\|\bar{\omega}_t\|_{b,2}\|\bar{u}_t\|_{b,4}\|\nabla\omega_t^1\|_{b,4} \leq
C\|\bar{\omega}_t\|_{b,2}\|\bar \omega_t\|_{b,k-1,2}\|\omega_t^1\|_{b,k,2}. 
\end{equation*}
We deduce that 
\begin{equation*}
\begin{aligned} d\|\bar{\omega}_t\|_{b,2}^2&+2\displaystyle\sum_{i=1}^{\infty}\langle
\bar{\omega}_t, \mathcal{L}_{i}\bar{\omega}_t \rangle_b dW_t^i
\leq C\|\omega_t^1\|_{b,k,2}\|\bar{\omega}_t\|_{b,2}^2dt. \end{aligned}
\end{equation*}
\noindent Similar arguments are used to control $\|\partial^\alpha \bar{\omega}_t\|_{b,2}^2$  where $\alpha$ is a multi-index with $|\alpha|\le k-1$ and to deduce that there exists a constant $C=C(R)$ such that 
\[
\begin{aligned} d\|\partial^\alpha\bar{\omega}_t\|_{b,2}^2&+2\displaystyle\sum_{i=1}^{\infty}\langle\partial^\alpha
\bar{\omega}_t, \partial^\alpha\mathcal{L}_{i}\bar{\omega}_t \rangle_b dW_t^i\leq
C\|\omega_t^1\|_{b,k,2}\|\bar{\omega}_t\|_{b,2}^2dt, \end{aligned}
\] 
where we use the control (see Lemma \ref{apriori})
\[
\big\langle \partial^\alpha\bar{\omega}_t, \partial^\alpha\mathcal{L}_{i}^2
\bar{\omega}_t \big\rangle_b  + \langle\partial^\alpha\mathcal{L}_{i} \bar{\omega}_t, \partial^\alpha\mathcal{L}_{i} \bar{\omega}_t \rangle_b\le C\|\bar\omega\|_{b,k,2}^2.
\]
We need to pay special attention to the case when $|\alpha|= k-1$, as  $\partial^\alpha\mathcal{L}_{i}^2
\bar{\omega}_t$ is no longer well defined. In this case we use the weak form of equation \eqref{truncatedEulere} to rewrite $\big\langle \partial^\alpha\bar{\omega}_t, \partial^\alpha\mathcal{L}_{i}^2\bar{\omega}_t \rangle_b$ as $-\big\langle \partial^{\alpha_1} \partial^\alpha\bar{\omega}_t, \partial^{\alpha_2}\mathcal{L}_{i}^2\bar{\omega}_t \rangle_b$ and then we can proceed as above by using that 
\[
-\big\langle \partial^{\alpha_1} \partial^\alpha\bar{\omega}_t, \partial^{\alpha_2}\mathcal{L}_{i}^2\bar{\omega}_t \rangle_b  + \langle\partial^\alpha\mathcal{L}_{i} \bar{\omega}_t, \partial^\alpha\mathcal{L}_{i} \bar{\omega}_t \rangle_b\le C\|\bar{\omega}_t\|_{b,k,2}^2.
\]
The above control is true for functions in $\mathcal{W}_b^{k+1,2}(\mathbb{T}^2)$ and, by the continuity of both sides in the above inequality, it is also true for functions which belong to the larger space $\mathcal{W}_b^{k,2}(\mathbb{T}^2)$, since $\mathcal{W}_b^{k+1,2}(\mathbb{T}^2)$ is dense in $\mathcal{W}_b^{k,2}(\mathbb{T}^2)$. The proof is now concluded in an identical manner as that for the uniqueness of the original equation.
    
\subsection{Existence of solution for the truncated equation}\label{etruncated}
Given the fact that the topology of $\mathcal{W}_b^{k,2}(\mathbb{T}^2)$ is equivalent to the topology of $\mathcal{W}^{k,2}(\mathbb{T}^2)$ for a nondegenerate weight $b$ (see \cite{KufnerWeightedSobSpaces}), the proof is similar to the proof from \cite{CL1}.
The strategy is to construct an approximating sequence of processes that will converge in distribution to a solution of the equation \eqref{truncatedEulere}. This justifies the existence of a weak solution. Together with the pathwise uniqueness obtained in Section \ref{utruncated}, we then deduce that a strong unique solution exists.  
For any $t\geq 0$ we construct the sequence $(\omega_t^{\nu_n,R,n})_{n\geq 0}$ with 
$\omega_t^{\nu_0,R,0} := \omega_0^0$ and for $n \geq 1$: 
\begin{equation}\label{itsyst} 
\begin{aligned} 
&\omega_0^{\nu_n,R,n} := \omega_0^n \\
& d{\omega}_t^{\nu_n,R,n}= \left(\nu_n\Delta \omega_t^{\nu_n,R,n} - f_R(u_t^{\nu_{n-1},R,n-1})\mathcal{L}_{u_t^{\nu_{n-1},R,n-1}} \omega_t^{\nu_n,R,n}\right)dt - \displaystyle %
\sum_{i=1}^{\infty}\mathcal{L}_i\omega_t^{\nu_n,R,n}\circ dW_t^{i,n},
\end{aligned} 
\end{equation}
where $\nu_n = \frac{1}{n}$ is the viscous parameter and  $u_t^{\nu_{n-1},R,n-1}=K\omega_t^{\nu_{n-1},R,n-1}.$ 
The corresponding It\^{o} form of equation \eqref{itsyst} is \footnote{%
The stochastic It\^{o} integral is understood here in the usual sense, see \cite{DaPratoZabczyk} .} 
\begin{equation}\label{itsystito}
d\omega_t^{\nu_n,R,n} = \left(\nu_n\Delta\omega_t^{\nu_n,R,n}+P_{t}^{n-1,n}(\omega_t^{\nu_n,R,n})\right)dt-\sum_{i=1}^{\infty} \mathcal{L}_{i}\omega_t^{\nu_n, R,n}dW_t^{i,n},   
\end{equation}
where $P_t^{n-1,n}(\omega_t^{\nu_n,R,n})$ is defined as  
\begin{equation}\label{ptqt}
P_t^{n-1,n}(\omega_t^{\nu_n,R,n}):= -f_R(u_t^{\nu_{n-1},R,n-1})\mathcal{L}_{u_t^{\nu_{n-1},R,n-1}}\omega_t^{\nu_n,R,n}+\frac{1}{2}\displaystyle\sum_{i=1}^{\infty}\mathcal{L}_{i}^2\omega_t^{\nu_n,R,n}, \ \ t\ge 0.
\end{equation}
\begin{theorem}\label{psns}If $\omega_0^{\nu_n,R,n} \in C^{\infty}(\mathbb{T}^2)$ then the two-dimensional stochastic vorticity equation \eqref{itsystito}
admits a unique  global $\mathcal F_t$-adapted solution
$\omega^{\nu_n,R,n}=\{\omega_t^{\nu_n,R,n},t\in [0,\infty)\}$ which belongs to the space 
$C\big([  0,\infty);C^\infty(\mathbb{T}^2) \big)$.  
\end{theorem}

\noindent The stochastic equation \eqref{itsystito} is a particular case of the more general equation $(1.1)-(1.2)$ analysed in Chapter 4, Section 4.1, pp.129 from \cite{Rozovskii}. All the assumptions required by Theorem 1 and Theorem 2 in \cite{Rozovskii}, Chapter 4, are fulfilled. 
Therefore there exists a unique solution $\omega_t^{\nu_n,R,n}$ which belongs to the class $L^2(0,T; \mathcal{W}_b^{k,2}(\mathbb{T}^2)) \cap C([0,T], \mathcal{W}_b^{k-1,2}(\mathbb{T}^2))$ 
and satisfies equation (\ref{itsystito}) for all $t\in[0,T]$ and for all $\omega$ in $\Omega' \subset \Omega$ with $\mathbb{P}(\Omega') = 1.$ 
Furthermore, since the conditions are fulfilled for all $k\in\mathbb{N}$, using Corollary 3
from pp. 141 in \cite{Rozovskii}, we obtain that $\omega_t^{\nu_n,R,n}$ is $\mathbb{P}$-a.s. in $C\big([0,T], C^{\infty}(\mathbb{T}^2)\big)$. 
Note that $u_t^{\nu_{n-1},R,n-1} \in C^{\infty}(\mathbb{T}^2)$ for any $n \geq 1$, using the regularity properties of the vorticity to stream function operator and an inductive argument. One has $u_t^{\nu_{n-1},R,n-1} = K \omega_t^{\nu_{n-1},R,n-1}$ where $K$ is the vorticity to stream function operator as before. We have 
\begin{equation*}
\|K\omega_t^{\nu_{n-1,R,n-1}}\|_{b, k+1, p} \leq C \|\omega_t^{\nu_{n-1,R,n-1}}\|_{b,k,p}.
\end{equation*}
Since $\omega_t^{\nu_{n-1}}$ belongs to $C^{\infty}(\mathbb{T}^2)$ we deduce that $u_t^{\nu_{n-1},R,n-1}$ also belongs to $C^{\infty}(\mathbb{T}^2)$. This, together with the initial assumptions \eqref{xiassumpt} ensures that all the coefficients of equation \eqref{itsystito} are infinitely differentiable. The uniform boundedness is ensured by the truncation $f_R(u_t^{\nu_{n-1},R,n-1})$, as proved in Lemma \ref{apriori}. 

\begin{remark}\label{cont} (\textit{Continuity of the approximating sequence}). 
There exists a constant $C=C(T)$ independent of $n$ and $R$ such that    
\begin{equation*}
 \mathbb{E}[\|\omega_t^{\nu_n,R,n} -
\omega_s^{\nu_n,R,n}\|_{L_b^2}^4] \leq C(t-s)^2,\ \ t,s\in[0,T].
\end{equation*}
In particular, by the Kolmogorov-\v{C}entsov
criterion (see \cite{KaratzasShreve}), the processes $\omega^{\nu_n,R,n}$ have continuous trajectories in $L_b^2(\mathbb{T}^2)$. An explicit proof of this fact is based on $L^2$-estimates for each term and can be found in \cite{CL1} Section 5.2.  
\end{remark} 

\begin{proposition}\label{tightnessI} 
The laws of the family of solutions $(\omega^{\nu_n,R,n})_{\nu_n\in [0,1]}$ is relatively compact in the space of probability measures over $D([0,T],L_b^{2}(\mathbb{T}^2))$ for any $T\geq 0$. 
\end{proposition}

\vspace{3mm}
\noindent \textbf{\textit{Proof of existence for the solution of equation \eqref{truncatedEulere} }} \\

  From  Proposition \ref{tightnessI} and the fact that $\displaystyle\lim_{n\rightarrow \infty} \omega_0^{\nu_n,R,n} 
=\omega_0$ we can deduce, using a diagonal subsequence argument, the existence of a subsequence $(\omega^{\nu_{n_j}})_{j}$ with $\displaystyle\lim_{j\rightarrow \infty}\nu_{n_j}=0$ 
which is convergent in distribution over  $D([0,\infty),L_b^2(\mathbb{T}^2))$. We show that the limit of the corresponding distributions is the distribution of a stochastic process which solves the truncated equation \eqref{truncatedEulere}. This justifies the existence of a weak probabilistic solution. By using the Skorokhod representation theorem (see \cite{Billingsley}), there exists a probability space $(\tilde{\Omega}, \tilde{\mathcal{F}}, \tilde{\mathbb{P}})$ and a sequence of processes $(\tilde{\omega}^{\nu_n,R,n}, \tilde{u}^{\nu_n,R,n}, (\widetilde{W}^{i,n})_i, n=1, \infty)$ which has the same distribution as that of the original converging subsequence and which converges, for $n \rightarrow \infty$, almost surely to a triplet $(\tilde{\omega}^R, \tilde{u}^R, (\widetilde{W}^{i,n})_i )$ in $D([0,T],L_b^2(\mathbb{T}^2)\times \mathcal{W}_b^{1,2}(\mathbb{T}^2)\times \mathbb{R}^{\mathbb N})$. Observe that $\omega^{\nu_n,R,n}$ and $\tilde{\omega}^{\nu_n,R,n}$ have the same distribution, therefore for any test function $\varphi \in C^{\infty}(\mathbb{T}^2)$ 
we have 
\begin{equation}\label{reff}
\begin{aligned}
\langle \tilde{\omega}_t^{\nu_n,R,n}, \varphi\rangle_b &= \langle \tilde{\omega}_0^{\nu_n,R,n}, \varphi\rangle_b + \nu_n\displaystyle\int_0^t \langle\tilde{\omega}_s^{\nu_n,R,n},  \Delta\varphi\rangle_b ds + \displaystyle\int_0^t f_R(\tilde{u}_s^{\nu_{n-1,R,n-1}})\langle \tilde{\omega}_s^{\nu_n,R,n},\mathcal{L}_{u_s^{\nu_{n-1},R,n-1}} \varphi \rangle_b ds \\
& + \frac{1}{2}\displaystyle\sum_{i=1}^{\infty}\displaystyle\int_0^t \langle\tilde{\omega}_s^{\nu_n,R,n},  \mathcal{L}_i^2\varphi\rangle_b ds + \displaystyle\sum_{i=1}^{\infty}\displaystyle\int_0^t \langle \tilde{\omega}_s^{\nu_n,R,n},  \mathcal{L}_{i}\varphi\rangle_b d\tilde{W}_s^{i,n}.
\end{aligned}
\end{equation}
Observe also that there exists a constant $C=C(R,T)$
such that  
\begin{equation}\label{C(R,T)}
\displaystyle\sup_{n\geq 1}\tilde{\mathbb{E}}\bigg[\displaystyle\sup_{s\in[0,T]}\|\tilde{\omega}^{\nu_n,R,n}_s\|_{b,k,2}^4\bigg] \leq C,
\end{equation}
where $\tilde{\mathbb{E}}$ is the expectation with respect to $\tilde{\mathbb{P}}$.  We prove this in Lemma \ref{apriori} for the original sequence, but since $\tilde{\omega}^{\nu_n,R,n}$ satisfies the same SPDE, the same a priori estimates hold for $\tilde{\omega}^{\nu_n,R,n}$. We know that the space of continuous functions is a subspace of the space of c\`adl\`ag functions, therefore the Skorokhod topology relativised to the space of continuous functions coincides with the uniform topology. It follows that the sequence  $(\tilde{\omega}^{\nu_n,R,n}, \tilde{u}^{\nu_n,R,n}, (\widetilde{W}^{i})_i, n=1, \infty)$ converges $\tilde{\mathbb{P}}$-almost surely to $(\tilde{\omega}^R, \tilde{u}^R, (\widetilde{W}^{i,n})_i )$ when $n \rightarrow \infty$, in the uniform norm. It also holds that
\[ 
\lim_{n\rightarrow\infty}\tilde{\mathbb{E}}\left[\int_0^t || \tilde{\omega}_s^{\nu_n,R,n}-\tilde{\omega}_s^R||^2ds\right]=0
\]
and since
\begin{equation}\label{tt}
\begin{aligned} 
\displaystyle\sum_{i=1}^{\infty}\displaystyle\tilde{\mathbb{E}}\left[\int_0^t \langle \tilde{\omega}_s^{\nu_n,R,n}-\tilde{\omega}_s^R, \mathcal{L}_i\varphi\rangle_b^2ds\right] &\le  \displaystyle\sum_{i=1}^{\infty} \|\mathcal{L}_{i}\varphi\|_{b,2}^2 \tilde{\mathbb{E}}\left[\int_0^t || \tilde{\omega}_s^{\nu_n,R,n}-\tilde{\omega}_s^R||^2ds\right] \\
&\leq C\|\varphi\|_{b,1,2}^2 \tilde{\mathbb{E}}\left[\int_0^t || \tilde{\omega}_s^{\nu_n,R,n}-\tilde{\omega}_s^R||^2ds\right]
\end{aligned} 
\end{equation}
also the limit of the right hand side of \eqref{tt} converges to $0$ (we use here the control 
$\sum_{i=1}^{\infty} \|\mathcal{L}_{i}\varphi\|_{b,2}^2 \leq C\|\varphi\|_{b,1,2}^2$ assumed in  \eqref{xiassumpt}).  
Now Theorem 4.2 in \cite{KurtzProtterII} allows us to conclude that the stochastic term 
$$
\sum_{i=1}^{\infty}\displaystyle\int_0^t \langle \tilde{\omega}_s^{\nu_n,R,n}, \mathcal{L}_i\varphi\rangle_b d\tilde W_s^{i,n}
$$
converges in distribution to 
$
\displaystyle\sum_{i=1}^{\infty}\displaystyle\int_0^t \langle \tilde{\omega}_s^R, \mathcal{L}_i\varphi\rangle_b d\tilde W_s^{i}.
$
Using a similar application of the Skorokhod representation theorem we can also assume that on the probability space $(\tilde{\Omega}, \tilde{\mathcal{F}}, \tilde{\mathbb{P}})$, the term   \\$\displaystyle\sum_{i=1}^{\infty}\displaystyle\int_0^t \langle \tilde{\omega}_s^{\nu_n,R,n},  \mathcal{L}_{i}\varphi\rangle_b d\tilde W_s^{i,n}$ converges to $\displaystyle\sum_{i=1}^{\infty}\displaystyle\int_0^t \langle \tilde{\omega}_s^R,  \mathcal{L}_{i}\varphi\rangle_b d\tilde W_s^i$
$\tilde{\mathbb{P}}$-almost surely (and also in $L_b^2(\tilde{\mathbb{P}})$).
The convergence of the remaining terms is shown in a similar fashion and similar to \cite{CL1}. The only difference is that we have weighted scalar products and norms, but given the fact that $b$ is a bounded nondegenerate function, all the convergences hold as requested. 

We have proven so far that  there exists a weak/distributional solution in the sense of  Definition \ref{solutions}. part b. on the probability space $(\tilde{\Omega}, \tilde{\mathcal{F}}, \tilde{\mathbb{P}})$. Moreover, since $\tilde{\omega}^R$ belongs to the space $\mathcal{W}_b^{k,2}(\mathbb{T}^2) \hookrightarrow C^{k-m}(\mathbb{T}^2)$ the solution is also strong, as a solution on the space $(\tilde{\Omega}, \tilde{\mathcal{F}}, \tilde{\mathbb{P}})$. It follows that  $(\tilde{\omega}, \tilde{u}, (\widetilde{W}^i)_i)$ is a weak/probabilistic solution of the truncated equation \eqref{truncatedEulere}  in the sense of Definition \ref{solutions} part c. Together with the pathwise uniqueness proved in Section \label{uniquenessEulert} \ref{uniquenessEulert} and using the Yamada-Watanabe theorem for the infinite-dimensional setting (see, for instance, \cite{Rockner}) we conclude the existence of a strong solution of the truncated great lake equation. Continuity follows as an application of the Kolmogorov-Čentsov criterion for the approximating process, which provides a control for $\mathbb{E}\left[\|\omega_t^{\nu_n,R,n}-\omega_s^{\nu_n,R,n}\|_{b,k,2}^4\right]$. Then we can pass to the limit using an argument similar to the one from Section \ref{contEuler}, to control $\mathbb{E}\left[\|\omega_t^{R}-\omega_s^{R}\|_{b,k,2}^4\right]$. Therefore, we have obtained a solution in the sense of  Definition \ref{solutions}, part a. Now using the embedding $\mathcal{W}_b^{k,2}(\mathbb{T}^2) \hookrightarrow C^{k-m}(\mathbb{T}^2)$ with $\ 2\leq m \leq k$ and $k\geq 4$ we conclude that the solution is classical when $k \geq 4$.

\section{Continuity with respect to initial conditions  \label{contEuler}}
In this section we prove Theorem \ref{mainresult2}.
Let $\omega$, $\tilde\omega$ be two $C( [  0,\infty);\mathcal{W}_b^{k,2}(\mathbb{T}^2) )  $-solutions of equation \eqref{maineqnito} and define the process $B=(B_t)_t$ such that  $B_t:=\displaystyle\int_0^t\|\omega_s\|_{b,k,2}ds$, for any $ t\ge 0$. Let $\omega^R$, $\tilde\omega^R$ be their corresponding truncated versions and also let $(\omega_t^{\nu_n,R,n})_{n\geq 0}$ and  $(\tilde\omega_t^{\nu_n,R,n})_{n\geq 0}$ be, respectively, the corresponding sequences constructed as in Section \ref{etruncated} on the same space after the application of the Skorokhod representation theorem. By Fatou's lemma, applied twice, we deduce that   
\begin{equation*} 
\begin{aligned}
\mathbb E [e^{-CB_t}||\omega_t-\tilde\omega_t||_{b,k-1,2}^2] & \le 
\mathbb E \big[\displaystyle\liminf_n e^{-CB^n_t}||\omega_t^{\nu_n,R,n}-\tilde\omega_t^{\nu_n,R,n}||_{b,k-1,2}^2\big] \\
& \leq  \displaystyle\liminf_n \mathbb{E} [e^{-CB^n_t}||\omega_t^{\nu_n,R,n}-\tilde\omega_t^{\nu_n,R,n}||_{b,k-1,2}^2]
\end{aligned}
\end{equation*}
where $B^n=(B_t^n)_t$ is the process defined by $B_t^n:=\displaystyle\int_0^t\|\omega_s^{\nu_n,R,n}\|_{b,k-1,2}ds$, for any $ t\ge 0$.
Following a similar proof with that of the uniqueness of the vorticity equation presented in Section \ref{uoriginal}, we deduce that  
there exists a positive constant $C$ independent of the two solutions and  independent of $n$ such that 
\[
E [e^{-CB^n_t}||\omega_t^{\nu_n,R,n}-\tilde\omega_t^{\nu_n,R,n}||_{b,k-1,2}^2]\le
||\omega_0-\tilde\omega_0||_{b,k-1,2}^2
\]
which gives the result. We emphasize that we use here the fact that the processes $(\omega_t^{\nu_n,R,n})_{n\geq 0}$ and $(\tilde\omega_t^{\nu_n,R,n})_{n\geq 0}$ take values in $\mathcal{W}_b^{k+2,2}(\mathbb{T}^2)$ as an essential ingredient, a property that was not true for either the solution of the original great lake equation \eqref{maineqnito} or its truncated version.  

\section{Relative compactness} \label{sectiontightness}
In this section we prove that the approximating sequence of solutions constructed in Section \ref{etruncated} is relatively compact in the space $D([0,\infty), L_b^2(\mathbb{T}^2))$. For this, we will use Kurtz' criterion for relative compactness (see \cite{EthierKurtz} Theorem 8.6). 
We denote by $(S^n(t))_t$ the semigroup of the generator $A:=\nu_n\Delta$.
This semigroup is strongly continuous (see \cite{Lunardi}) and for any $f \in L^2(\mathbb{T}^2)$
it is true that 
\begin{equation*}
\|S^n(t)f\|_{k,2} \leq \|f\|_{k,2}.
\end{equation*}

\noindent \textbf{\textit{Proof of Proposition \ref{tightnessI}}} \\
  We need to show that, for every $\eta>0$, there exists a compact set $K_{\eta, t} \subset L_b^2(\mathbb{T}^2)$ such
that $\displaystyle\sup_{n}\mathbb{P}\big(\omega_t^{\nu_n,R,n}\notin
K_{\eta, t}\big) \leq \eta.$ We define the compact
\[
K_{\eta, t} :=\left\{ \omega \in\mathcal{W}_b^{k,2}(\mathbb{T}^2) | \ \  \|\omega\|_{b,k,2} <\left(\frac{C}{\eta}\right)^{\frac{1}{4}} \right\}
\] 
where $C$ is the constant which appears in the a priori estimates  (Lemma \ref{apriori}).
 By a Sobolev compact embedding theorem,  $K_{\eta, t}$ is a compact set in
$L_b^2(\mathbb{T}^2)$ and 
\[
\displaystyle\sup_{n}\mathbb{P}\big(\omega_t^{\nu_n,R,n}\notin
K_{\eta, t}\big) = \displaystyle\sup_{n}\mathbb{P}\left( \|\omega_t^{\nu_n,R,n}\|_{b,k,2} \geq \left(\frac{C}{\eta}\right)^{\frac{1}{4}} \right)\le\displaystyle\sup_{n}\frac{\eta}{C}\mathbb{E}\left[\sup_{t\in [0,T]}\|\omega_t^{\nu_n,R,n}\|_{b,k,2}^4
\right]\le\eta.  
\]
In order to show relative compactness, we need to justify part b)\ of Kurtz' criterion, as per Theorem 8.6 in \cite{EthierKurtz}. We will  show that there exists a  family $(\gamma_{\delta}^n)_{0<\delta<1}$ of nonnegative random
variables such that 
\[
\mathbb{E} \big[\|\omega_{t+l}^{\nu_n,R,n} - \omega_t^{\nu_n,R,n}\|_{b,2}^2   | \mathcal{F}_t^{}\big] \leq \mathbb{E}\big[%
\gamma_{\delta}^{n} | \mathcal{F}_t^{} \big]
\] 
and $\displaystyle \lim_{\delta \rightarrow 0} \displaystyle\lim
\sup_{n} \mathbb{E}\big[ \gamma_{\delta}^{n} \big] =0 $ for $t \in
[0, T]$. We use the mild form of equation (\ref{itsystito}), that is
\begin{equation*}
\omega_t^{\nu_n,R,n}=S^n(t)\omega_0^{\nu_n,R,n} -
\displaystyle\int_{0}^{t}S^n(t-s)P_{s}^{n-1,n}(\omega_s^{\nu_n,R,n}))ds-\sum_{i=1}^{\infty} \int_{0}^{t}S^n(t-s)\mathcal{L}_{i}\omega_s^{\nu_n, R,n}dW_s^{i,n},
\end{equation*}
with $P_s^{n-1,n}$ as defined in (\ref{ptqt}) and $S^n(t):=e^{\nu_n\Delta t}$.
One has 
\begin{equation*}
\begin{aligned} &\|\omega_{t+l}^{\nu_n,R,n} - \omega_t^{\nu_n,R,n}\|_{b,2}^2  \leq
C\bigg( \|(S^n(t+l)-S^n(t))\omega_0^{\nu_n,R,n}\|_{b,2}^2 \\ & +
\bigg\|\displaystyle\int_0^t(S^n(t+l-s)-S^n(t-s))P_s^{n-1,n}(\omega_s^{\nu_n,R,n})ds\bigg\|_{b,2}^2 +
\bigg\|\displaystyle\int_t^{t+l}S^n(t+l-s)P_s^{n-1,n}(\omega_s^{\nu_n,R,n})ds\bigg\|_{b,2}^2 \\ & +
\bigg\|\displaystyle\sum_{i=1}^{\infty} \int_{0}^{t}(S^n(t+l-s)-S^n(t-s))\mathcal{L}_{i}\omega_s^{\nu_n, R,n}dW_s^{i,n}\bigg\|_{b,2}^2 +
\bigg\|\displaystyle\sum_{i=1}^{\infty} \int_{t}^{t+l}(S^n(t+l-s)\mathcal{L}_{i}\omega_s^{\nu_n, R,n}dW_s^{i,n}\bigg\|_{b,2}^2 \bigg)
\end{aligned}
\end{equation*}

\noindent The calculations follow as per \cite{CL1} Section 6. 
The process $\gamma_l^{\nu_n}$ is defined as 
\begin{equation*}
\begin{aligned} \gamma_l^{\nu_n} &:= \|(S^n(l)-1)\omega_0^{\nu_n,R,n}\|_{b,2} +
\displaystyle\int_0^T \|(S^n(l)-1)P_s^{n-1,n}(\omega_s^{\nu_n,R,n})\|_{b,2}^2 ds+ Cl^2\displaystyle\sup_{s
\in [0, T+1]}\|P_s^{n-1,n}(\omega_s^{\nu_n,R,n})\|_{b,2}^2 \\ &+ l^2 \displaystyle\sup_{s \in[0, T+1]}
\|\mathcal{L}_i\omega_s^{\nu_n,R,n}\|_{b,2}^2 + \displaystyle\int_0^T
\|\big((t+l-r)^{\alpha-1}S^n(l)-(t-r)^{\alpha-1}\big)S^n(t-r)z(r)\|_{b,2}^2dr
\end{aligned}
\end{equation*}
and $\gamma_{\delta}^{n} := \displaystyle\sup_{l \in [0,
\delta]}\gamma_l^{\nu_n}. $ From Lemma \ref{apriori} we deduce that there exist two constants $c_1$ and $c_2$ such that $\mathbb{E}[%
\displaystyle\sup_s\|P_s^{n-1,n}(\omega_s^{\nu_n,R,n})\|_{b,2}^2] \leq c_1$
and $\mathbb{E}[\displaystyle\sup_s\|\mathcal{L}_i\omega_s^{\nu_n,R,n}\|_{b,2}^2] \leq c_2.$ The integrands in the integrals above converge
pointwise to $0$ when $l \rightarrow 0$ due to the strong continuity of the
semigroup $S^n$. 
At the same time, they are bounded by integrable functions,
therefore the convergence is uniform in space by the dominated convergence
theorem. Then the requirement 
\begin{equation*}
\displaystyle \lim_{\delta \rightarrow 0} \displaystyle \sup_{\nu_n} 
\mathbb{E}\big[ \gamma_{\delta}^{n} \big] =0
\end{equation*}
is met. In conclusion all the conditions required by Kurtz' criterion are
fulfilled and therefore $(\omega_t^{\nu_n, R, n})_{\nu_n}$ is relatively compact. 

\section{A priori estimates} \label{apriorisection}

\begin{lemma}\label{apriori}$\left.\right. $\\ 
Let $\omega_t$ be the solution of the vorticity equation \eqref{maineqnito} and $\omega_t^{\nu_n,R,n}$ the solution of the linear approximating equation \eqref{itsyst}. Then the following properties hold: 
\begin{enumerate}
\item[i.] For any $f\in \mathcal{W}_b^{2,2}(\mathbb{T}^2)$ we have
\[
\big\langle f, \mathcal{L}_{i}^2
f \big\rangle_b  + \langle\mathcal{L}_{i} f, \mathcal{L}_{i} f \rangle_b = 0.
\]
\item[ii.] If  $\omega_t, \omega_t^{\nu_n,R,n}\in \mathcal{W}_b^{k,2}(\mathbb{T}^2)$  then the following transport formulae hold $\mathbb{P}$-almost surely:
\begin{equation*}
\begin{aligned} 
&\|\omega_t\|_{L_b^p} =\|\omega_0\|_{L_b^p} \\
&\|\omega_t^{\nu_n,R,n}\|_{L_b^p}\leq
\|\omega_0\|_{L_b^p} \\
& \|\omega_t\|_{\infty} = \|\omega_0\|_{\infty} \\
& \|\omega_t^{\nu_n, R, n}\|_{\infty} \leq \|\omega_0\|_{\infty}. 
\end{aligned}
\end{equation*}
provided that the initial condition is in $L^{\infty}(\mathbb{T}^2)$. 
\item[iii.] If  $\omega\in \mathcal{W}_b^{k,2}(\mathbb{T}^2)$, then  $(P_t^{n-1,n})_t$ defined in \eqref{ptqt} and $(\mathcal{L}_i\omega_t^{\nu_n,R,n})_t$ are processes with paths taking values in $L_b^2(\mathbb{T}^2)$. 

\item[iv.] There exists a constant $C_1$ such that: 
\begin{equation*}
\big|\big\langle \partial^k \omega_t^{\nu_n, R, n}, \partial^k \big( \mathcal{L}_i^2\omega_t^{\nu_n, R, n}\big)\big\rangle_b + \big\langle\partial^k\big(\mathcal{L}_i\omega_t^{\nu_n, R, n}\big),\partial^k\big(\mathcal{L}_i\omega_t^{\nu_n, R, n}%
\big)\big\rangle_b \big| \leq C_1\|\omega_t^{\nu_n, R, n}\|_{b,k,2}^2.
\end{equation*}
\item[v.] There exist some constants $C_2$
and $C_2^{\prime}$ such that: 
\begin{equation*}
\big|\big\langle \partial^k \omega_t^{\nu_n, R, n}, f_R(u_t^{\nu_{n-1}, R, n-1}) \partial^k\big(\mathcal{L}_{u_t^{\nu_{n-1}, R, n-1}}
\omega_t^{\nu_n, R, n}\big)\big \rangle_b \big| \leq
C_2\|\partial^{k}u_t^{\nu_n, R, n-1}\|_{b,2}^a\|u_t^{\nu_n, R, n-1}\|_{b,2}^{1-a}\|\omega_t^{\nu_n, R, n}%
\|_{b,k,2}^2
\end{equation*}
with $0<a\leq 1$, and 
\begin{equation*}
|\langle \partial^k \omega_t^{\nu_n, R, n}, f_R(u_t^{\nu_{n-1}, R, n-1})\partial^k\big(\mathcal{L}_{u_t^{\nu_{n-1}, R, n-1}}
\omega_t^{\nu_n, R, n}\big) \rangle_b| \leq C_2^{\prime }\|\omega_t^{\nu_n, R, n}\|_{b,k,2}^2.
\end{equation*}%

\end{enumerate}
\end{lemma}
The technical details of the proof of the above inequalities are similar to the ones proven in Lemma 19 from \cite{CL1}. The only difference is given by the fact the on the right hand side we always have also the norm of the bounded nondegenerate function $b$, which does not alter qualitatively the result. However, a slightly more significant difference appears in the control of the higher derivatives of the Lie derivative $\mathcal{L}_{u_t^{\nu_{n-1},R,n-1}}\omega_t^{\nu_n,R,n}$ where instead of the divergence-free assumption from \cite{CL1} we use the regularity properties proven in Section \ref{biotsavartsection}.

\begin{proposition}\label{aprioriapprox}
There exists a constant $\mathcal{C}(R,T)$ independent of $n$ such that 
\begin{equation*}
\mathbb{E}\big[\displaystyle\sup_{t \in [0,T]}\|\omega_t^{\nu_n,R,n}\|_{b,k,2}^4\big] \leq \mathcal{C}(R,T). 
\end{equation*}
\end{proposition}

\noindent\textbf{Proof}:  
After applying the It\^{o} formula to the approximating sequence \eqref{itsyst} we obtain 
\begin{equation*}
\begin{aligned}
\|\partial^k\omega_t^{\nu_n, R, n}\|_{b,2}^2 &= \|\partial^k\omega_0^{\nu_n, R, n}\|_{b,2}^2 + 2\nu_n\displaystyle\int_0^t\langle \partial^k\omega_s^{\nu_n, R, n}, \partial^{k+2}\omega_s^{\nu_n, R, n}\rangle_b ds \\
&- 2\displaystyle\int_0^t \langle \partial^k\omega_s^{\nu_n, R, n}, f_R(u_s^{\nu_{n-1},R, n-1})\partial^k\mathcal{L}_{u_s^{\nu_{n-1}, R, n-1}}\omega_s^{\nu_n, R, n}\rangle_b ds \\
& + \displaystyle\sum_{i=1}^{\infty}\displaystyle\int_0^t \langle \partial^k\omega_s^{\nu_n,R,n}, \partial^k\mathcal{L}_i^2\omega_s^{\nu_n,R,n}\rangle_b ds \\ 
& + \displaystyle\sum_{i=1}^{\infty}\displaystyle\int_0^t\langle \partial^k\mathcal{L}_i\omega_s^{\nu_n, R, n}, \partial^k\mathcal{L}_i\omega_s^{\nu_n, R, n}\rangle_b ds \\
& - 2\displaystyle\sum_{i=1}^{\infty}\displaystyle\int_0^t\langle \partial^k\omega_s^{\nu_n, R, n}, \partial^k\mathcal{L}_i\omega_s^{\nu_n, R, n} \rangle_b dW_s^{i,n}. \\
\end{aligned}
\end{equation*}
We analyse each term. One can write
\begin{equation*}
\langle \partial^k\omega_s^{\nu_n,R,n},\partial^{k+2}\omega_s^{\nu_n,R,n}\rangle_b = -\|\partial^{k+1}\omega_s^{\nu_n,R,n}\|_{b,2}^2 \leq 0. 
\end{equation*}
We want to estimate the other terms independently of $\nu_n$. All terms are estimated above. After summing up we have
\begin{equation*}
\begin{aligned}
\mathbb{E}\big[\displaystyle\sup_{s\in[0,t]}\|\omega_s^{\nu_n, R, n}\|_{b,k,2}^2\big] &\leq \mathbb{E}\big[\|\partial^k\omega_0^{\nu, R, n}\|_{b,2}^2\big] + C_2'(T)\displaystyle\int_0^t\mathbb{E}\big[\displaystyle\sup_{s\in[0,t]}\|\omega_s^{\nu_n,R,n}\|_{b,k,2}^2\big] ds \\
& + C_1\displaystyle\int_0^t\mathbb{E}\big[\displaystyle\sup_{s\in[0,t]}\|\omega_s^{\nu_n,R,n}\|_{b,k,2}^2\big]ds
 + 2\mathbb{E}\bigg[\displaystyle\sup_{s\in[0,t]}\displaystyle\sum_{i=1}^{\infty}\displaystyle\int_0^t\langle \partial^k\omega_s^{\nu_n, R, n}, \partial^k\mathcal{L}_i\omega_s^{\nu_n, R, n} \rangle_b dW_s^{i,n}\bigg]. \\
\end{aligned}
\end{equation*}
Let 
\begin{equation*}
B_t := \displaystyle\sum_{i=1}^{\infty}\displaystyle\int_0^t\langle \partial^k\omega_s^{\nu_n, R, n}, \partial^k\mathcal{L}_i\omega_s^{\nu_n, R, n} \rangle_b dW_s^{i,n} \ \ \ \ \hbox{and} \ \ \ \
\beta_t := \|\omega_t^{\nu_n,R,n}\|_{b,k,2}^2.
\end{equation*}
$B_t$ is a local martingale. We have
\begin{equation*}
\beta_t \leq \beta_0 - 2B_t + (C_2'+C_1) \displaystyle\int_0^t\beta_sds 
\end{equation*}
and by Gronwall lemma
\begin{equation*}
\mathbb{E}\big[\displaystyle\sup_{s\in[0,t]}\beta_s^2\big] \leq 2e^{2(C_2'+C_1) t}\bigg(\beta_0^2 + 4\mathbb{E}\big[\displaystyle\sup_{s\in[0,t]}B_s^2\big]\bigg).
\end{equation*}
Using the Burkholder-Davis-Gundy inequality there exists a constant $\tilde{\alpha}$ such that
\begin{equation*}
\mathbb{E}\big[\displaystyle\sup_{s\in[0,t]}B_s^2\big] \leq \tilde{\alpha} \mathbb{E}\bigg[\sqrt{\langle B\rangle_s^2}\bigg].
\end{equation*}
It is easy to show that $\mathbb{E}\big[\langle B\rangle_s\big] \leq C_2^2\mathbb{E}[\displaystyle\sup_{s\in[0,t]} \beta_s^2]$: 
following the same calculations as those from step v. Lemma \ref{apriori} we can see that
\begin{equation*}
\langle \partial^k\omega_s^{\nu_n, R, n}, \partial^k\mathcal{L}_i\omega_s^{\nu_n, R, n} \rangle_b \leq C_2\|\omega_s^{\nu_n,R,n}\|_{b,k,2}^2
\end{equation*}
if instead of making use of the truncation corresponding to $\|u_s^{\nu_{n-1},R,n-1}\|_{b,k,2}$ we take into account assumption \eqref{xiassumpt}. Hence
\begin{equation*}
\langle B \rangle_s \leq C_2^2 \displaystyle\sup_{s\in[0,t]}\|\omega_s^{\nu_n,R,n}\|_{b,k,2}^4 = C_2^2 \displaystyle\sup_{s\in[0,t]}\beta_s^2
\end{equation*} 
and the above inequality follows. In conclusion 
\begin{equation*}
\mathbb{E}\big[\displaystyle\sup_{s\in[0,t]} \beta_s^2\big] \leq \tilde{C_1}(T)\beta_0^2 + \tilde{C_2}(T)\displaystyle\int_0^t \mathbb{E}\big[\displaystyle\sup_{s\in[0,t]} \beta_s^2]ds
\end{equation*}
with $\tilde{C}_1(T):= 2\exp(2(C_2^{'}+C_1)T)$ and $\tilde{C}_2(T):=8\tilde{\alpha}C_2^2\exp(2(C_2^{'}+C_1)T)$. Using again Gronwall's inequality we have
\begin{equation*}
\mathbb{E}\big[\displaystyle\sup_{s\in[0,T]}\|\omega_s^{\nu_n,R,n}\|_{b,k,2}^4\big] \leq \mathcal{C}(R,T)
\end{equation*}
with $\mathcal{C}(R,T):= \tilde{C}_1(T)\beta_0^2\exp(\tilde{C}_2(T)T)$. \\

\noindent The following lemma is essential when showing that the limit of the approximating sequence satisfies the great lake equation in $\mathcal{W}_b^{k,2}(\mathbb{T}^2)$ although the relative compactness property holds in $D\big([0,T], L_b^{2}(\mathbb{T}^2)\big)$ for all $T>0$.  
\noindent Given the fact that the proof is similar to the proof of Lemma 23 in \cite{CL1}, we do not redo it here. 
  
\begin{lemma}\label{lim}$\left.\right.$\\
i. Assume that $(a_n)_n$
is a sequence of functions such that $\displaystyle\lim_{n\rightarrow\infty} a_n=a$ in $L_b^2(\mathbb{T}^2)$ and $\displaystyle\sup_{n>1}
\left\Vert a_n\right\Vert_{b,s,2}<\infty$ for $s \geq 0$. Then $a \in \mathcal{W}_b^{s,2}(\mathbb{T}^2)$ and $\left\Vert a \right\Vert_{b,s,2}<\displaystyle\sup_{n>1}\left\Vert a_n\right\Vert_{b,s,2}$. Moreover, $\displaystyle\lim_{n\rightarrow\infty}a_n = a$ in $\mathcal{W}_b^{s',2}(\mathbb{T}^2)$ for any $s'<s$.\\
ii.  
Assume that $a_n:\Omega\rightarrow \mathcal{W}_b^{s,2}(\mathbb{T}^2)$ is a sequence of measurable maps such that,  $\displaystyle\lim_{n\rightarrow\infty} a_n=a$ in $L_b^2(\mathbb{T}^2)$, $\mathbb{P}$-almost surely or $\displaystyle\lim_{n\rightarrow\infty} a_n=a$ in distribution. Further assume that $\displaystyle\sup_{n>1}\mathbb E[\left\Vert a_n\right\Vert_{b,s,2}^{2}]<\infty$. Then, $\mathbb{P}$-almost surely, $a\in \mathcal{W}_b^{s,2}(\mathbb{T}^2) $ and $\mathbb E[\left\Vert a\right\Vert_{b,s,2}^{2}] \le \displaystyle\sup_{n>1}\mathbb E[\left\Vert a_n\right\Vert_{b,s,2}^{2}]$. \\
\end{lemma}

\begin{remark}\label{normequiv}
The norm $\|\cdot\|_{b,m,2}$ is equivalent to the norm defined as $||| f |||:= \|f\|_{b,2} + \|D^{m}f\|_{b,2}$, therefore it is enough to show that all properties hold for the $L_b^2$ norm of $f$ and for the $L_b^2$ norm of the maximal derivative $D^mf$ (see \cite{Brezis} pp. 217).  
\end{remark}

\vspace{3mm}
\noindent\textbf{Acknowledgements} \\

\noindent During this work, Dan Crisan has been partially supported by the ERC Synergy grant STUOD - DLV-856408. Oana Lang has been supported by the EPSRC grant EP/L016613/1, through the Mathematics of Planet Earth Centre for Doctoral Training, and the EPSRC grant EP/N023781/1.
The authors would like to thank Professor Darryl Holm for having suggested this problem, and also Erwin Luesink, Michael James Leahy, Colin Cotter, Peter Jan van Leeuwen, and Sebastian Reich for many constructive discussions held with them during the preparation of this work.

\end{document}